\documentstyle[12pt]{article}
\newcommand{\gl}{gl(M|N)}
\newcommand{\SY}{super Yangian\ }
\newcommand{\Y}{Y(\gl)}
\newcommand{\lam}{\lambda}
\newcommand{\Lam}{\Lambda}
\newcommand{\ba}{\begin{eqnarray}}
\newcommand{\na}{\end{eqnarray}}
\newcommand{\ban}{\begin{eqnarray*}}
\newcommand{\nan}{\end{eqnarray*}}
\newtheorem{lemma}{Lemma}

\newtheorem{theorem}{Theorem}

\newcommand{\x}{\otimes}

\begin{document}
\title{\small{\bf THE $gl(M|N)$ SUPER YANGIAN AND ITS \\
FINITE DIMENSIONAL REPRESENTATIONS }}
\author{\small R. B. ZHANG\\
\small  Department of Pure Mathematics\\
\small  University of Adelaide\\
\small  Adelaide, S. A.,  Australia}
\maketitle

\begin{abstract}
Methods are developed for systematically constructing the finite
dimensional irreducible representations of the super
Yangian $\Y$ associated with the Lie superalgebra $\gl$.
It is also shown that every finite
dimensional irreducible representation of $\Y$ is of highest weight type,
and is uniquely characterized by a highest weight.  The necessary and
sufficient conditions for an irrep to be finite dimensional are given. \\

\noindent
{\bf Mathematics subject classifications(1991)}: 17B37, 81R50, 17A70\\
{\bf Running title}: Representations of $\Y$
\end{abstract}

\normalsize

\vspace{3cm}
\noindent

\section{ Introduction}

This is the second of a series of papers developing the representation
theory of the super Yangians.
In an earlier publication\cite{I}, the finite dimensional irreducible
representations of the super Yangian $Y(gl(1|1))$ were classified, and
explicit bases for such representations were constructed.
It is the aim of this paper to carry out a similar program for the
super Yangian associated with the Lie superalgebra $\gl$ for all
$M$ and $N$.  In particular, we will prove that every finite dimensional
irrep of $\Y$ is of highest weight type, and  uniquely characterized
by a highest weight, and give the necessary and sufficient conditions
for an irrep to be finite dimensional. Methods will also be developed
for systematically constructing the finite dimensional irreps.

The structures and representations of the Yangians associated with
ordinary Lie algebras were studied extensively\cite{Drinfeld}
\cite{Tarasov}\cite{Cherednik}\cite{Kirillov}\cite{Chari}
\cite{Reshetikhin}\cite{Molev}.  The central result in the
representation theory is Drinfeld's Theorem\cite{Drinfeld},
which provides a characterization of finite dimensional irreps in terms
of highest weight polynomials, and also gives  the necessary and
sufficient conditions for an irrep to be finite dimensional.
However, more detailed structural information on Yangian irreps,
such as their dimensions and their decompositions with respect to
the underlying Lie algebras,  is difficult to obtain,
primarily due to the existence of finite dimensional
indecomposable but non -  irreducible representations.
At present, the understanding of the structures of irreps of Yangians
is still incomplete.

Super Yangians and their quantum analogues arose naturally from the
algebraic description of integrable lattice models with Lie
superalgebra symmetries. Their representation theory  plays a
central role in the study of such models, e.g., by applying the
algebraic Bethe Ansatz to diagonalize the transfer matrices.
Super (quantum)Yangians are also algebraic systems of considerable
mathematical interest;   their structure and representations merit a
thorough investigation in their own right.

Some structural features of super (quantum)Yangians, in particular,
their connections with the Lie superalgebras and the related quantum
supergroups,  have already been studied by Nazarov\cite{Nazarov} and
also in \cite{II}.  The finite dimensional irreps of the simplest
\SY $Y(gl(1|1))$ have been classified, and explicit bases for such
irreps have been constructed\cite{I}. It is the aim of this letter
to systematically develop the representation theory of the super
Yangian $\Y$.

In section $2$, we define the super Yangian $\Y$, then study its structure.
In particular, we will prove a $BPW$ theorem, which will be our starting
point for developing the representation theory.  Section $3$  contains
the main results. We prove that every finite dimensional irrep
of $\Y$ is of highest weight type, and is uniquely characterized
by the highest weight.  A general method is developed to construct highest
weight irreps; the necessary and sufficient conditions for an irrep to
be finite dimensional are also given. In section 4 we outline another
construction of the finite dimensional irreps of $\Y$, which is similar
to Kac' induced module construction\cite{Kac} for Lie superalgebras.
This construction should be useful for investigating detailed structures
of irreps.

This letter relies heavily on \cite{I}.  The reasoning and techniques
employed in proving some of the results of subsection 3.1
are adopted from that paper.   To make this letter reasonably
self contained, we spell out these proofs in some detail.

\section{ Structure of $\Y$ }
\subsection{ Definition }
The \SY    $\Y$ was first defined  by Nazarov\cite{Nazarov} following
the Faddeev - Reshetikhin - Takhatajan formalism of quantization of
algebraic structures. We present the definition here to fix our notation.

We will work on the complex number field ${\bf C}$. The underlying
vector space of the Lie superalgebra $\gl$ is ${\bf Z}_2$  graded,
with a homogeneous basis  $\{ E^a_b |$ $ a, b = 1, 2, ..., M+N\}$.
Introduce the gradation index $[\ ]:  \{1, 2, ..., M+N\}\rightarrow
{\bf Z}_2$ such that
$[a]=\left\{\begin{array}{l r}
          0,& a\le M,\\
          1,& a>M.
           \end{array} \right. $
Let $\gl_\theta$, $\theta\in {\bf Z}_2$,  be the vector space over ${\bf C}$
spanned by the $E^a_b$ with $[a]+[b]\equiv$ $\theta \ (mod\ 2)$.  Then $\gl_0$
and $\gl_1$ are the even and odd subspaces of $\gl$ respectively. We will
abuse the notation a bit and define $[\ ]: \gl_0\cup \gl_1$  $\rightarrow$
${\bf Z}_2$, $[x]=\left\{\begin{array}{l r}
          0,& x\in\gl_0,\\
          1,& x\in\gl_1.
           \end{array} \right. $
The Lie superalgebra $\gl$ is this ${\bf Z}_2$ graded vector space endowed
with the bilinear  graded bracket
$[\ , \ \}: \gl\otimes\gl\rightarrow\gl$,
\ba
{[} E^a_b, \ E^c_d\} &=& \delta_b^c E^a_d -
(-1)^{([a]+[b])([c]+[d])}  E^c_b\delta_d^a.
\na
For convenience, we can regard $\gl$ as embedded in its universal
enveloping algebra. Thus the graded bracket $[, \ \}$ can be interpreted
as the graded commutator
\ba
[x, \ y\}&=&x y - (-1)^{[x][y]} y x.\label{commutator}
\na

The vector  module of $\gl$ is an $(M+N)$ - dimensional
${\bf Z}_2$-graded vector space $V$, spanned by the homogeneous elements
$\{ v^a | a = 1, 2, ..., M+N\}$ where $v^a$ is even if $[a]=0$ and odd
otherwise. The action of $\gl$  on $V$ is defined by $E^{a}_{b} v^c =
\delta_b^c v^a$. We denote the associated vector representation of
$\gl$ by $\pi$. Then in this basis $\pi(E^{a}_{b})=e^{a}_{b}$, where
$e^{a}_{b}\in End(V)$ are the standard matrix units.

Define the permutation operator $P: V\otimes V\rightarrow V\x V$ by
$$P(v^a\x v^b)=(-1)^{[a][b]} v^b\x v^a.$$
Then explicitly, we have
\ban
P&=&\sum_{a, b=1}^{M+N}e^a_b \otimes e^b_a (-1)^{[b]}.
\nan
It is well known that the following $R$ matrix
\ba
R(u)&=& 1 + { {P}\over {u}},\   \ \ \ \ u\in{\bf C},
\na
satisfies the graded Yang - Baxter equation. \\

Let us introduce
\ban
L(u)&=&\sum_{a, b=1}^{M+N}(-1)^{[b]}t^a_b(u)\otimes e^b_a,\\
t^a_b(u)&=&(-1)^{[b]}\delta^a_b + \sum_{n=1}^{\infty} t^a_b[n] u^{-n},
\nan
where $u$ is an indeterminate, and the
$t^a_b[n]$, $0<n\in{\bf Z}_+$, are homogeneous elements of $\Y$
such that $t^a_b[n]$ is even if $[a]+[b]\equiv 0(mod\ 2)$ and odd otherwise.
The $L(u)$ belongs to the ${\bf Z}_2$ graded vector space $\Y\otimes
End(V)$ and is even.  Now  $\Y$ is the ${\bf Z}_2$ graded associative algebra
generated by the $t^a_b[n]$, $0<n\in{\bf Z}_+$, with the following
defining relations
\ba
L_1(u) L_2(v) R_{12}(v-u) &=& R_{12}(v-u)  L_2(v) L_1(u) .\label{def}
\na
Note that equation (\ref{def}) lives in $\Y\otimes End(V)\otimes End(V)$.
The multiplication of the factors on both sides are defined with respect
to the grading of this triple tensor product.
To gain some concrete feel about this algebra, we put (\ref{def}) into a more
explicit form:
\ban
{[} t^{a_1}_{b_1}(u), t^{a_2}_{b_2}(v)\} &=&
{ {(-1)^{\eta(a_1, b_1; a_2, b_2)}}\over {u-v} }
\left[ t^{a_2}_{b_1}(u)t^{a_1}_{b_2}(v)
- t^{a_2}_{b_1}(v) t^{a_1}_{b_2}(u)\right],  \\
{\eta}(a_1, b_1; a_2, b_2)&\equiv&
[a_1][a_2] + [b_1]([a_1]+[a_2]) (mod \ 2);
\nan
or equivalently,
\ba
{[} t^{a_1}_{b_1}[m], t^{a_2}_{b_2}[n]\}=
\delta^{a_2}_{b_1}t^{a_1}_{b_2}[m+n-1] - (-1)^{([a_1]+[b_1])([a_2]+[b_2])}
t^{a_2}_{b_1}[m+n-1]\delta^{a_1}_{b_2}\nonumber\\
+ (-1)^{\eta(a_1, b_1; a_2, b_2)} \sum_{r=1}^{min(m, n)-1}
\left\{ t^{a_2}_{b_1}[r]t^{a_1}_{b_2}[m+n-1-r]
-t^{a_2}_{b_1}[m+n-1-r]t^{a_1}_{b_2}[r]\right\} .\label{modes}
\na

$\Y$  admits co - algebraic structures compatible with
the associative multiplication.  We have the
co - unit $\epsilon: \Y\rightarrow{\bf C}$,
$t^a_b[k]\mapsto\delta_{0 k}\delta^a_b(-1)^{[a]}$,
the co - multiplication $\Delta: \Y\rightarrow \Y\otimes \Y$,
$L(u)\mapsto L(u)\otimes L(u)$, and also
the antipode $S: \Y\rightarrow\Y$, $L(u)\mapsto L^{-1}(u)$.
Thus $\Y$ is indeed a ${\bf Z}_2$ graded Hopf algebra.

Note that $\Y$ also admits the following generalized tensor product structure:
\ba
\Delta^{(k-1)}_{\alpha}:  \Y&\rightarrow&\Y^{\otimes k}, \nonumber \\
L(u)&\mapsto& L(u+\alpha_1)\otimes L(u+\alpha_2)\otimes ...
\otimes L(u+\alpha_k), \label{comultiplication}
\na
where  $\alpha_1 =0$, and $\alpha_i$, $i=2, 3,  ..., k$,  are
a set of arbitrary complex parameters.  Explicitly, we have
\ban
\Delta^{(k-1)}_{\alpha}(t^a_b(u))&=&\sum_{a_1, ..., a_{k-1}}
(-1)^{\sum_{i=1}^{k-1}\{[a_i]+([a_0]+[a_i])([a_i]+[a_{i+1}])\} }\\
&\times& t^{a_1}_b(u)\otimes t^{a_2}_{a_1}(u+\alpha_2)\otimes ...
\otimes t^a_{a_{k-1}}(u+\alpha_k)
\nan
where $a_0=b$, and $a_k=a$.

Another  useful fact is the existence of an
automorphism $\phi_f: \Y$ $\rightarrow$ $\Y$ associated with
each power series $f(x)$ $=$ $1 + f_1 x^{-1} + f_2 x^{-2}+...$,
which is defined by
\ba
t^a_b(x)&\mapsto&{\tilde t}^a_b(x)=f(x)t^a_b(x).    \label{auto}
\na
As can be easily seen, the ${\tilde t}^a_b$ satisfy exactly the same
relations as the $t^a_b$ themselves. \\

Some further simple properties of the \SY are worth observing.
Note that $\Y$ as a Hopf algebra is a deformation\cite{Deform}
of the universal enveloping algebra of the infinite dimensional Lie
superalgebra $\widehat\gl^{(+)}$  $=$  $\gl\otimes {\bf C}[[t]]$,
where ${\bf C}[[t]]$ denotes the ring of polynomials in the
indeterminate $t$.
Set $E^a_b[k]=E^a_b\otimes t^k$, $k\in{\bf Z}_+$.
Then the graded bracket for $\widehat\gl^{(+)}$ reads
\ba
{[} E^a_b[k], \ E^c_d[l]\} &=& \delta_b^c E^a_d[k+l] -
(-1)^{([a]+[b])([c]+[d])}  E^c_b[k+l]\delta_d^a. \label{loop}
\na
The universal enveloping algebra  $U(\widehat\gl^{(+)})$ of this Lie
superalgebra is a ${\bf Z}_2$-graded associative algebra, which may
be thought as generated by $E^a_b[k]$,  $k\in{\bf Z}_+$, subject to
the relations (\ref{loop}) but with the left hand side interpreted as
the graded commutator defined in (\ref{commutator}).

This algebra in fact has the structure of a ${\bf Z}_2$-graded Hopf algebra.
In particular, its co - multiplication is given by
\ban
\delta: U(\widehat\gl^{(+)})&\rightarrow&
U(\widehat\gl^{(+)})\otimes U(\widehat\gl^{(+)}),\\
E^a_b[k]&\mapsto& E^a_b[k]\otimes 1 + 1\otimes  E^a_b[k].
\nan

In order to see that $\Y$ is indeed a deformation of the Hopf superalgebra
$U(\widehat\gl^{(+)})$, we set $t^{a}_{b}[m+1]=\kappa^{-m}E^{a}_{b}[m]$,
$m\in{\bf Z}_+$, where $\kappa$ is an indeterminate.  Then $\Y$ is
isomorphic to the algebra $\tilde U$  generated by $E^{a}_{b}[m]$, $m\in{\bf
Z}_+$,
subject to the relations
\ban
& &{[} E^{a_1}_{b_1}[m], E^{a_2}_{b_2}[n]\}=
\delta^{a_2}_{b_1}E^{a_1}_{b_2}[m+n] - (-1)^{([a_1]+[b_1])([a_2]+[b_2])}
E^{a_2}_{b_1}[m+n]\delta^{a_1}_{b_2}\nonumber\\
&+\kappa &(-1)^{\eta(a_1, b_1; a_2, b_2)} \sum_{r=1}^{Min(m, n)}
\left\{ E^{a_2}_{b_1}[r-1]E^{a_1}_{b_2}[m+n-r]
-E^{a_2}_{b_1}[m+n-r]E^{a_1}_{b_2}[r-1]\right\}.
\nan
Regard $\tilde U$ as an algebra defined on the polynomial ring
${\bf C}[[\kappa]]$. Then it is clear from the above equation
that $U(\widehat\gl^{(+)})=\tilde U/\kappa\tilde U$. Also,
the co - multiplication $\Delta$ of $\Y$ induces
a co - associative co - multiplication $\Delta: \tilde U
\rightarrow  \tilde U$, which is clearly the deformation of
the co - multiplication $\delta$ of $U(\widehat\gl^{(+)})$.

Important structural and representation theoretical properties of $\Y$
can be obtained by investigating this Hopf superalgebra within the
framework of deformation theory\cite{Deform}. Results will be reported
in a future publication. \\

The \SY  $\Y$ contains several subalgebras, which will be useful
for developing the representation theory.  We can easily see
from the defining relations (\ref{modes}) that the generators
$t^a_b[1]$ form the Lie superalgebra $\gl$. For each $n>1$,
the generators $t^a_b[n]$ transform as the components
of an adjoint tensor operator of this $\gl$.
Define a map $t^a_b[n]\mapsto\delta_{1 n}E^a_b$,
where $E^a_b$ are the standard generators of $\gl$.
Then it extends to an algebra homomorphism $\Y\rightarrow U(\gl)$.

There also exist various Hopf (super) subalgebras of $\Y$. In particular,
the following will be used in the remainder of the letter:
\ba
Y(gl(M)) & generated\ by & \{t^a_b[n] |\  a, b=1, 2, ..., M,
\ 0<n\in{\bf Z}_+\}; \nonumber\\
Y(gl(N)) & generated\ by &  \{t^a_b[n] | a, b=M+1, M+2, ..., M+N,
\ 0<n\in{\bf Z}_+\};   \nonumber\\
Y(gl(1|1))& generated\ by &
\{t^M_{M+1}[n], \ t^{M+1}_M[n], \  t^M_M[n], \ t^{M+1}_{M+1}[n]
| \ 0<n\in{\bf Z}_+\}. \label{subalgebras}
\na
Note that although both $Y(gl(M))$ and $Y(gl(N))$ are even
subalgebras, together they do not form a subalgebra of $\Y$.
This leads to certain complications in the development of
the representation theory.

\subsection{BPW theorem}
We now prove a version of the BPW theorem for the \SY $\Y$, which will
be of crucial importance for developing the representation theory.
Let us introduce a filtration on $\Y$.  Define
the degree of a generator $t^a_b[n]$ by $deg(t^a_b[n])=n$, and
require that the degree of a monomial
$t^{a_1 }_{b_1}[n_{1 }]t^{a_2}_{b_2}[n_{2 }]
...t^{a_k}_{b_k}[n_{k}]$ is $\sum_{r=1}^k n_{r}$.
Let $Y_p$ be the vector space over $\bf C$ spanned by monomials of
degree not greater than $p$. Then
\ban
...\supset Y_p\supset Y_{p-1}\supset ... \supset Y_1\supset Y_0=\bf C,\\
Y_p Y_q\subset Y_{p+q}.
\nan

Let $z_1, \ z_2, ..., \ z_k$  be some $t^a_b[n]$'s. Consider the product
$Z=z_1  z_2 ... z_k$, which is assumed to have $deg(Z)=p$.  It directly
follows from the defining relations (\ref{modes}) of the \SY  $\Y$  that
for any permutation $\sigma $ of $(1, 2, ..., k)$,
\ban
z_1  z_2 ... z_k -
\epsilon(\sigma) z_{\sigma(1)} z_{\sigma(2)} ... z_{\sigma(k)}
\nan
belongs to $Y_{p-1}$, where $\epsilon(\sigma)$ is
$-1$ if $\sigma$ permutes
the odd elements in $z_1, \ z_2, ..., \ z_k$ an odd number of times,
and $+1$ otherwise.  In particular, if $t^a_b[n]$ is odd, then
$(t^a_b[n])^2\in Y_{2n-1}$.
Therefore,  given any ordering of the
generators $t^a_b[n]$,  $0<n\in{\bf Z}_+$,   $ a, b\in\{1, 2, ..., M+N\}$,
their ordered products of degrees less or equal to $p$ span $Y_p$,
where the products do not contain factors $(t^a_b[n])^2$ if
$[a]+[b]\equiv 1 (mod \ 2)$.  It immediately follows that the ordered
products of all degrees span the underlying vector space of $\Y$.
As we will show presently, the ordered products are also linearly
independent, thus form  a basis for  $\Y$.

Define $U_p=Y_p/Y_{p-1}$. Then the multiplication of $\Y$ defines
a bilinear map $U_p\otimes U_q\rightarrow U_{p+q}$, which
extends to a multiplication $U\otimes U\rightarrow U$ for the
space $U=\oplus_{p=0}^{\infty}U_p$, turning $U$ into an associative algebra.
This algebra is isomorphic to the algebra of polynomials
$G[X]$ in the variables $X^a_b[n]$, $a,  b\in \{1,2, ..., M+N\}$,
$n=1, 2, ...$, with the isomorphism $U\cong G[X]$ defined by
\ban
t^a_b[n]&\mapsto& X^a_b[n], \ \ \ \  \forall a, \ b, \ n,
\nan
where $X^a_b[n]$ is an ordinary indeterminate if $[a]+[b]\equiv 0 (mod \ 2)$,
and is a Grassmannian variable if $[a]+[b]\equiv 1 (mod \ 2)$.
(Note that for any Grassmannian variables $\zeta_1$ and $\zeta_2$
we have $\zeta_i\zeta_j = - \zeta_j \zeta_i$, $i, j =1, 2$. )
Since monomials in $X^a_b[n]$ are linearly independent as elements of $G[X]$,
we conclude that ordered products of the $t^a_b[n]$
(not allowing powers of order higher than $1$ of the odd $t^a_b[n]$)
are linearly independent.  To summarize, we have
\begin{theorem}:
For any given ordering of the generators $t^a_b[n]$,
$ a, b\in\{1, 2, ..., M+N\}$, $0<n\in{\bf Z}_+$,
the ordered products of the $t^a_b[n]$ containing no second and higher order
powers of the odd generators  form a basis of $\Y$.
\end{theorem}

It is useful to construct an explicit  basis for $\Y$.  To do that, we need
to fix  some notations.   Consider the pairs $(a,  b)$,
$a,  b\in \{1,2, ..., M+N\}$.
Let
\ban
\Phi_+&=&\{(a, b)| a<b\},\\
\Phi_-&=&\{(a, b)| a>b\},\\
\Phi_0&=&\{(a, a)\},\\
\Phi_\pm^{(\theta)}&=&\{ (a, b)\in\Phi_\pm| [a] +[b] \equiv \theta (mod \ 2)\}.
\nan
Given any $p=(a, b)$, we denote $\bar p=(b, a)$.
We introduce  a total ordering $\succ$ $(\prec)$ of all the pairs
in the following way: for any $p_+\in \Phi_+$, $p_0\in \Phi_0$,
$p_-\in \Phi_-$, we define $p_+\succ p_0 \succ p_-$, or equivalently,
$p_-\prec p_0 \prec p_+ $. For $(a, b), (c, d) \in \Phi_+$, we define
$(a, b)\succ (c, d)$ (i.e., $(c, d)\prec (a, b)$) if $a<c$ or $a=c, \ b>d$;
for  $p_1,  p_2$ belonging to $\Phi_-$, $p_1\succ p_2$ if
$\bar{p_1}\prec\bar{p_2}$;
and $(a, a)\succ(b, b)$ if $a<b$.

Let
\ban
Q^{\{k\}}_{(a, b)}[\{n\}]&=&
(t^a_b[n_1])^{k_1}(t^a_b[n_2])^{k_2}...(t^a_b[n_r])^{k_r},\ \ \
n_1<n_2<...<n_r,\\
k_1, ..., k_r&\in&\left\{\begin{array}{l l}
                   {\bf Z}_+, & [a] +[b] \equiv 0 (mod \ 2),\\
                   \{0, 1\}, & [a] +[b] \equiv 1 (mod \ 2).
                   \end{array}\right.
\nan
Define the ordered product
$$\prod_{p}^\succ Q^{\{k_{p}\}}_{p}[\{n_{p}\}],$$
which positions $Q^{\{k_{p}\}}_{p}[\{n_{p}\}]$ on the right
of $Q^{\{k_{p'}\}}_{p'}[\{n_{p'}\}]$ if $p\succ p'$.
Now
\begin{lemma}
The following elements
\ba
\prod_{p\in\Phi_-^{(1)}}^\succ Q^{\{k_{p}\}}_{p}[\{n_{p}\}]
 \prod_{q\in\Phi_-^{(0)}}^\succ Q^{\{k_{q}\}}_{q}[\{n_{q}\}]
 \prod_{r\in\Phi_0}^\succ Q^{\{k_r\}}_r[\{n_r\}]
 \prod_{s\in\Phi_+^{(0)}}^\succ Q^{\{k_s\}}_s[\{n_s\}]
 \prod_{t\in\Phi_+^{(1)}}^\succ Q^{\{k_t\}}_t[\{n_t\}],\label{basis}
\na
form a basis of $\Y$.
\end{lemma}

\section{Finite Dimensional  Irreps}
We study structures of the finite dimensional irreps of $\Y$ in
this section.  Some of the results reported here are generalizations
of those on the representations of $Y(gl(1|1))$\cite{I} to the present case,
and the proofs of these results are also adopted from \cite{I}.

\subsection{Highest weight irreps}
Let $V$ be  an irreducible $\Y$ - module.
A nonzero element $v^\Lam_+\in V$ is called maximal if
\ba
t^a_b[n] v^\Lam_+=0, &\forall (a, b)\in\Phi_+, &n>0, \nonumber\\
t^a_a[n] v^\Lam_+=\lam_a[n] v^\Lam_+,
&a=1, 2, ..., M+N, & n>0,\label{high}
\na
where $\lam_a[n]\in{\bf C}$.  An irreducible module is called
a highest weight module if it admits a maximal vector.  We define
\ba
\Lam(x)=(\lam_1(x), \lam_2(x), ..., \lam_{M+N}(x)),
& \lam_a(x)=(-1)^{[a]}+\sum_{k>0}\lam_a[k]x^{-k}, \label{Lam}
\na
and call $\Lam(x)$ a highest weight of $V$.

Note that commutators amongst the $t^a_a[n]$ do not close on these
generators themselves; the usual practice of using  Lie's Theorem
to prove the existence of a common eigenvector for them does not
work here.   But nevertheless, we have
\begin{theorem}\label{existence}
Every finite dimensional irreducible $\Y$ - module $V$ contains a
unique(up to scalar multiples) maximal vector $v^\Lam_+$.
\end{theorem}

In order to prove the theorem, we need the following
\begin{lemma}\label{technical}
Define $V_0=\{v\in V\ |  \ t^a_b[n] v=0,  \forall (a, b)\in\Phi_+,
\  n>0\}$. Then \\
1). $V_0\ne 0$;\\
2). the generators $t^a_a[n]$  stabilize $V_0$; \\
3). for all $v\in V_0$,
\ban
[t^a_a[m], \ t^b_b[n]]v=0, &  \forall  a, b,  \  m, n.
\nan
\end{lemma}
{\em Proof}:
1). Let $E$ be an $(M+N)$ - dimensional vector space over ${\bf C}$.
We define a map $\alpha: \Phi_+\cup\Phi_-\rightarrow E$
by setting $\alpha(a, b)$ to  the vector with the $a-th$ component
being $+1$ and the $b - th$ component being $-1$, e.g.,
$\alpha(1, 3)=(1, 0, -1, 0,...,0)$.
We can introduce a partial ordering of all the elements of $E$ by
requiring that $\mu\succ\nu$ (i.e., $\nu\prec\mu$) if
$\mu=\nu+\sum_{p\in\Phi_+} l_p\alpha(p)$, $l_p\in{\bf Z}_+$,
where at least for one $p\in\Phi_+$, $l_p>0$.

If the nonvanishing element $v\in V$ is a  common eigenvector of the
$t^a_a[1]$,   then
\ban
t^a_a[1]v=\mu_a v, & \mu_a\in{\bf C}.
\nan
In this way,  every such $v$ is associated with a unique vector $\mu=$
$(\mu_1, \mu_2..., \mu_{M+N})\in E$, which we call the $\gl$ weight of $v$.
It is obvious that if a set of common eigenvectors all have different
$\gl$ weights, then they must be linearly independent.
Since $\{t^a_a[1] | a=1, 2, ..., M+N\}$ forms an abelian Lie algebra,
it follows Lie's Theorem that there exists at least one nonvanishing
element of $V$, which is their common eigenvector.

Let $v\in V$ be a common eigenvector of the $t^a_a[1]$ with a $\gl$ weight
$\mu$.  If $v\in V_0$, we have proved the first part of the Lemma.
If $v\not\in V_0$, by applying $t^a_b[n]$, $(a, b)\in \Phi_+$,  to $v$
we will arrive at other  common eigenvectors of the $t^a_a[1]$, which
have $\gl$ weights $\succ\mu$. Since $V$ is finite dimensional, there can
only exist a finite number of vectors with $\gl$ weights $\succ\mu$. Hence
repeated applications of the generators $t^a_b[n]$, $(a, b)\in \Phi_+$,
to $v$ will lead to a nonvanishing $v_0\in V$ such that
\ban
t^a_b[n]v_0&=&0, \ \ \  \forall (a, b)\in\Phi_+,\ \  \  n>0, \\
t^a_a[1]v_0&=&\lam_a[1] v_0, \ \ \   \forall a.
\nan
This proves that $V_0$ contains at least one nonzero element.

2). Let $v$ be a vector of $V_0$. We want to prove that
all $t^{a_k}_{a_k}[n_k]t^{a_{k-1}}_{a_{k-1}}[n_{k-1}]...t^{a_1}_{a_1}[n_1] v$,
$a_i=1, 2, ..., M+N$, $n_i>0$, $k\ge 0$,  are annihilated by
$t^a_b[m]$, $(a, b)\in\Phi_+$, $m>0$.
The $k=0$ case requires no proof. Assume that all the vectors
$v_{l}=t^{a_{l}}_{a_{l}}[n_{l}]...t^{a_1}_{a_1}[n_1] v$, $l<k$,
are in $V_0$. Then
\ban
&i).& (a, b)\in \Phi_+, \  \ \ b>c, \\
& &t^a_b[m]t^c_c[n_k]v_{k-1}=-[t^c_c[n_k], \ t^a_b[m]]v_{k-1}\\
& &=(-1)^{[c]+1}\sum_{r=0}^{min(m, n_k)-1} \left( t^a_c[r] t^c_b[m+n_k-1-r]
 -t^a_c[m+n_k-1-r]t^c_b[r]\right)v_{k-1}\\
& &=0,\\
&ii).& (a, b)\in \Phi_+, \  \ \ b\le c, \\
& & t^a_b[m]t^c_c[n_k] v_{k-1} = [t^a_b[m], \ t^c_c[n_k]]v_{k-1}\\
& &=\sum_{r=0}^{min(m, n_k)-1}
\left(t^c_b[r] t^a_c[m+n_k-1-r]
   -t^c_b[m+n_k-1-r]t^a_c[r]\right)v_{k-1}\\
& & \times    (-1)^{[a][c]+[b]([a]+[c])} =0.
\nan

3). The following defining relations of $\Y$
\ban
[t^a_a[m], \ t^a_a[n]]&=&0, \ \ \ \ \  \ \ a=1, 2, ..., M+N, \\
{[}t^b_b[m],\ t^c_c[n]]&=&(-1)^{[b]}\sum_{r=0}^{min(m, n)-1}
\left(t^c_b[r]t^b_c[m+n-1-r] - t^c_b[m+n-1-r]t^b_c[r]\right), \\
&  &  (b, c)\in\Phi_+,
\nan
and part $2)$ of the Lemma directly lead to part $3)$. \\

\noindent {\em Proof of theorem \ref{existence}}:
By part $3)$ of  Lemma \ref{technical},   the action of
the $t^a_a[n]$ on $V_0$ coincides with an abelian subalgebra of $gl(V_0)$.
Therefore, Lie's Theorem can be applied, and we conclude that
there exists at least one common eigenvector of all the
$t^a_a[n]$ in $V_0$. This proves the existence
of the highest weight vector. Assume that $v_+$ and $v_+'$ are two
highest weight vectors of $V$, which are not proportional to each other.
Applying $\Y$ to them generates two nonzero submodules of $V$,
which are not equal. This contradicts the irreducibility of $V$. \\

We now turn to the construction of highest weight irreps of $\Y$.
Let $N^+$, $N^-$ and $Y^0$ be the vector spaces spanned by the ordered
products (as defined in (\ref{basis})) of the elements $t^a_b[n]$,
with $(a, b)\in\Phi^+ $, $(a, b)\in\Phi^-$ and $(a, b)\in\Phi^0$
respectively.  Set $Y^+=Y^0N^+$, and $Y^-=N^-Y^0$.
We emphasize that these vector spaces do not form subalgebras of $\Y$.

Consider a one dimensional vector space ${\bf C}v^\Lam_+$.
We define a linear action of $Y^+$ on it by
\ba
t^a_b[n] v^\Lam_+&=&0,\ \ \ \  (a, b)\in\Phi_+, \nonumber\\
t^a_a[n] v^\Lam_+&=&\lam_a[n] v^\Lam_+,\nonumber\\
t^a_a[n]y_0 v^\Lam_+&=&\lam_a[n]y_0 v^\Lam_+, \ \ \ \ \forall y_0\in Y^0.
\label{highest}
\na
{}From the proof of  Lemma \ref{technical}  we can see that the definition
(\ref{highest}) is consistent with the commutation relations of $\Y$.
Now we define  the following  vector space
\ban
{\bar V}(\Lam)&=&\Y\otimes_{Y^+}v^\Lam_+.
\nan
Then ${\bar V}(\Lam)$ is a $\Y$ module, which is obviously isomorphic
to $N^-\otimes v^\Lam_+$.

The action of $\Y$ on this module is defined in the following way.
Every vector of ${\bar V}(\Lam)$ can be expressed
as $y\otimes v^\Lam_+$ for some  $y\in N^-$. For simplicity,
we write it as $y v^\Lam_+$.
Given any $u\in\Y$, $u y$ can be expressed as a linear sum of the
basis elements (\ref{basis}).  We write
\ban
u y&=& \sum a_{\alpha, \beta} y_-^{(\alpha)} y_0^{(\beta)}
+\sum b_{\mu, \nu, \sigma} y_-^{(\mu)}y_0^{(\nu)}y_+^{(\sigma)},\\
&\ & y_-^{(\alpha)}, y_-^{(\mu)}\in N^-,
\  y_0^{(\beta)}, y_0^{(\nu)} \in Y^0,
\  \ 1\ne y_+^{(\sigma)}\in N^+.
\nan
Then
\ban
u (y v^\Lam_+)&=& \sum a_{\alpha, \beta} y_-^{(\alpha)}
y_0^{(\beta)} v^\Lam_+.
\nan

The $\Y$ module ${\bar V}(\Lam)$ is infinite dimensional.
Standard arguments show that it is indecomposable, and
contains a unique maximal proper submodule $M(\Lam)$.  Construct
\ban
V(\Lam)&=& {\bar V}(\Lam)/M(\Lam).
\nan
Then $V(\Lam)$ is an irreducible highest weight $\Y$ module.  \\

Let $V_1(\Lam)$ and $V_2(\Lam)$ be two irreducible $\Y$ modules with
the same highest weight $\Lam(x)$.  Denote by $v^\Lam_{1,+}$ and
$v^\Lam_{2,+}$ their maximal vectors respectively.
Set $W=V_1(\Lam) \oplus V_2(\Lam)$. Then
$v^\Lam_+=(v^\Lam_{1,+}, v^\Lam_{2,+})$
is maximal, and repeated applications of $\Y$ to $v^\Lam_+$ generate
an $\Y$ submodule $V(\Lam)$ of $W$.  Define the module homomorphisms
$P_i: V(\Lam)\rightarrow V_i(\Lam)$  by
\ban
P_1(v_1, v_2)&=&(v_1, 0),\\
P_2(v_1, v_2)&=&(0, v_2), \  \  v_1\in V_1(\Lam), \ v_2\in V_2(\Lam).
\nan
Since
\ban
P_1(v^\Lam_{1,+}, v^\Lam_{2,+})&=&(v^\Lam_{1,+}, 0),\\
P_2(v^\Lam_{1,+}, v^\Lam_{2,+})&=&(0, v^\Lam_{2,+}),
\nan
it follows the irreducibility of $V_1(\Lam)$ and $V_2(\Lam)$ that
$Im P_i = V_i(\Lam)$.  Now $Ker P_1$ is a submodule of $V_2(\Lam)$.
The irreducibility of $V_2(\Lam)$ forces  either $Ker P_1=0$ or
$Ker P_1=V_2(\Lam)$.  But the latter case is not possible, as
$(0, v^\Lam_{2,+})\not\in W$. Similarly we can
show that $Ker P_2=0$.   Hence, $P_i$ are $\Y$ module isomorphisms.

To summarize the preceding discussions, we have
\begin{theorem}\label{uniqueness}
Corresponding to each  $\Lam(x)$ of the form (\ref{Lam}),
there exists a unique irreducible highest weight $\Y$ module
$V(\Lam)$ with highest weight $\Lam(x)$.  \\
\end{theorem}

Before closing this subsection, we consider some useful facts about
tensor products of  highest weight irreps of $\Y$.
Let $W(\mu)$ and $W(\nu)$ be finite dimensional irreducible $\Y$ modules
with highest weights $\mu(x)=(\mu_1(x), \mu_2(x), ..., \mu_{M+N}(x))$
and $\nu(x)$ $=$ $(\nu_1(x), \nu_2(x),  ..., \nu_{M+N}(x))$ respectively.
(The existence of such modules will be proven in the next subsection.)
Let $w^\mu_+$ and $w^\nu_+$ be the maximal vectors of these modules.
Then $v_+=w^\mu_+\otimes w^\nu_+$
$\in$ $W(\mu)\otimes W(\nu)$ is a maximal vector such that
\ban
t^a_a(x) v_+ &=& (-1)^{[a]}\mu_a(x)\nu_a(x) v_+,  \ \forall a.
\nan
Set $\lam_a(x)=(-1)^{[a]}\mu_a(x)\nu_a(x)$, and define the $\star$ -
 product of the highest weight $\mu(x)$ and $\nu(x)$ by
\ban
\mu(x)\star\nu(x)&=&(\lam_1(x), \lam_2(x), ..., \lam_{M+N}(x)).
\nan
Applying $\Y$ to $v_+$ generates an indecomposable $\Y$ module
${\bar V}(\mu\star\nu)$, which is contained  in $W(\mu)\otimes W(\nu)$.
The quotient module of ${\bar V}(\mu\star\nu)$ by its unique
maximal invariant submodule yields an irreducible $\Y$ module
$V(\mu\star\nu)$ with highest weight $\mu(x)\star\nu(x)$,
which is necessarily finite dimensional.
(The maximal invariant submodule of ${\bar V}(\mu\star\nu)$ can  of course
be trivial. In that case, ${\bar V}(\mu\star\nu)$ $=$ $V(\mu\star\nu)$.)
Clearly this discussion can be generalized to tensor products of more
than two irreps.

\subsection{Finite dimensionality conditions}
Let  $V(\Lam)$ be a finite dimensional irreducible
$\Y$ module with  highest weight $\Lam(x)$.
Denote its maximal vector by $v^\Lam_+$.
We now consider the actions of the subalgebras of $\Y$ defined by
(\ref{subalgebras}) on the maximal vector of $V(\Lam)$.
The following $V(\Lam)$  subspaces
$Y(gl(M))v^\Lam_+$,
$Y(gl(N))v^\Lam_+$  and
$Y(gl(1|1))v^\Lam_+$
furnish indecomposable modules over the subalgebras $Y(gl(M))$, $Y(gl(N))$,
and $Y(gl(1|1))$ respectively. It is obvious but very important to note that
these modules, being subspaces of $V(\Lam)$,  are finite dimensional.
Thus, the necessity part of the following theorem immediately
follows from Drinfeld's Theorem\cite{Drinfeld} (also see\cite{Molev})
 and a result of \cite{I}:
\begin{theorem}\label{iff}
The irreducible highest weight $\Y$ - module $V(\Lam)$
is finite dimensional if and only if its  highest weight
$\Lam(x)$ satisfies the following conditions
\ba
{ {\lam_a(x)}\over{\lam_{a+1}(x)} }&=&{ {P_a(x+(-1)^{[a]})}\over{P_a(x)}},
\ \ \ \ 1\le a < N+M,  \  a\ne M,  \nonumber\\
{ {\lam_M(x)}\over{\lam_{M+1}(x)} }&=&{ {\tilde Q_M(x)}\over{Q_M(x)}},
\label{result}
\na
where
\ban
P_a(x)=\prod_{i=1}^{K_a}(x+p_a^{(i)}), & 1\le a < N+M, \   a\ne M,\\
\tilde Q_M(x)=\prod_{i=1}^{K_M}\left( 1+ {{r_1^{(i)}}\over{x}}\right),\\
Q_M(x)=-\prod_{i=1}^{K_M}\left( 1+ {{r_2^{(i)}}\over{x}}\right),
&  \tilde Q_M(x), \ Q_M(x)     \   co - prime.
\nan
\end{theorem}
{\em Proof}: We only need to prove sufficiency.
Let us consider the special case where the highest weight
$\mu(x)=(\mu_1(x),\mu_2(x)..., \mu_{M+N}(x))$  of a $\Y$ irrep is of the form
\ba
\mu_a(x)&=&(-1)^{[a]}+\mu_a x^{-1},  \ \ \ \forall a.\label{Lie}
\na
We denote the irreducible $\Y$ module with highest weight $\mu(x)$
by $W(\mu)$, and the associated irrep by $\pi_\mu$.
This irreducible representation can be explicitly constructed using  the
irreducible $\gl$ representation $\gamma_\mu$
with highest weight $\mu=(\mu_1, \mu_2, ..., \mu_{M+N})$.
We have
\ban
\pi_\mu(t^a_b(u))&=&(-1)^{[b]}\delta^a_b + \gamma_\mu(E^a_b) u^{-1},
\ \ \forall a, b,
\nan
where $E^a_b$ are the standard $\gl$ generators.

In this case, the given conditions of the theorem are equivalent to
the following constraints on $\mu$:
\ba
\mu_a -\mu_{a+1}\in{\bf Z}_+, & 1\le a < N+M, \ a\ne M.  \label{lexical}
\na
This is nothing else but the necessary and sufficient conditions in order for
the $\gl$ irrep $\gamma_\mu$ to be finite dimensional.  Therefore, the
$\Y$ irrep $\pi_\mu$ is also finite dimensional, and we have proved the
sufficiency in this case. This also proves the fact that finite dimensional
irreps of $\Y$ indeed exist.

The next step in the proof  involves showing that up to
an automorphism $\phi_f$ of $\Y$
defined by (\ref{auto}),  every finite dimensional irreducible
$\Y$ module $V(\Lam)$ with highest weight $\Lam(x)$ sits inside
the tensor product of some irreducible $\Y$ modules
$W(\mu)$, where the highest weights of these modules are of the
form (\ref{Lie}). To do this, we note that  an  automorphism
$\phi_f$ transforms the highest weight
according to $\Lam(x)\mapsto f(x)\Lam(x)$.
Thus, it leaves the $P_a$,  $Q_M$, and $\tilde Q_M$ intact,
but an appropriate choice of $f(x)$
will change the components of $\Lam(x)$ into
polynomials in $x^{-1}$ defined by
\ban
\lam_a(x)&=&\prod_{c=1}^{a-1}Q_c(x) \prod_{d=a}^{M+N-1}\tilde Q_d(x), \ \
\forall a,
\nan
where we have used the following  notation
\ban
Q_a(x)=(-1)^{(K_a+1)[a]}\prod_{i=1}^{K_a}q_a^{(i)}(x),&
q_a^{(i)}(x)=(-1)^{[a]}\left(1+{{p_a^{(i)}}\over{x}}\right),\\
\tilde Q_a(x)=(-1)^{(K_a+1)[a]}\prod_{i=1}^{K_a}\tilde q_a^{(i)}(x),
&\tilde q_a^{(i)}(x)=(-1)^{[a]}\left(1+{{p_a^{(i)}+(-1)^{[a]}}
\over{x}}\right),  \ \  a\ne M, \\
\tilde q_M^{(i)}(x)=\left(1+{{r_1^{(i)}}\over{x}}\right),&
q_M^{(i)}(x)=-\left(1+{{r_2^{(i)}}\over{x}}\right).
\nan
Define
\ban
\mu^{(t, i)}(x)=(\mu^{(t, i)}_1(x),  \mu^{(t, i)}_2(x), ...,
\mu^{(t, i)}_{M+N}(x)), \\
\mu^{(t, i)}_a(x)=\left\{\begin{array}{l l}
                   q_t^{(i)}(x), & t<a, \\
                   \tilde q_t^{(i)}(x), & t\ge a,
                  \end{array}\right.\\
i=1, 2, ..., K_t, \  t=1, 2, ..., M+N-1.
\nan
Then $\Lam(x)$ can be expressed as the $\star$ -  product of all the
$\mu^{(t, i)}(x)$
\ban
\Lam(x)&=&\star_{t=1}^{M+N-1}\star_{i=1}^{K_t}\mu^{(t, i)}(x).
\nan

The $\mu^{(t, i)}(x)$ are of the form (\ref{Lie}) and satisfy the
conditions (\ref{lexical}).
Therefore, corresponding to each $\mu^{(t, i)}(x)$,
there exists a finite dimensional irreducible $\Y$ module $W(\mu^{(t, i)})$
with highest weight $\mu^{(t, i)}(x)$. It follows  the discussions at the
end of the last subsection that the tensor product module
$\otimes_{t=1}^{M+N-1}\otimes_{i=1}^{K_t}W(\mu^{(t, i)})$ (The order
of the tensor product is not important for us.)
contains as a submodule an indecomposable $\Y$ module $\bar V(\Lam)$.
A quotient module of $\bar V(\Lam)$ gives rise to an
irreducible $\Y$ module  which is isomorphic to $V(\Lam)$.
Being a submodule of a finite dimensional module, $\bar V(\Lam)$
is finite dimensional, and so is also $V(\Lam)$. This proves the
sufficiency of the given conditions of the theorem.

\section{ Another Construction of Irreps}
Experiences with the representation theories of the Lie superalgebras
and quantum supergroups urge us to ask whether a method similar
to Kac' induced module construction\cite{Kac} can be
developed to build irreps of $\Y$. The answer to this question is
affirmative. We now outline the construction.

Introduce an auxiliary algebra $Y(gl(M))\dot+Y(gl(N))$, which
is the product of the two Yangians $Y(gl(M))$ and $Y(gl(N))$.
In more explicit terms,  $Y(gl(M))\dot+Y(gl(N))$ is generated by
$\{\phi^i_j[n], \ \psi^\mu_\nu[n] |$ $i,   j=1, 2, ..., M$,
$\mu, \nu=M+1, M+2, ..., M+N$, $0<n\in{\bf Z}_+\}$,
where the $\phi^i_j[n]$ satisfy the standard defining relations of
a $gl(M)$ Yangian algebra $Y(gl(M))$, and the $\psi^\mu_\nu[n]$
satisfy relations of $Y(gl(N))$, while
\ban
[\phi^i_j[m],  \ \psi^\mu_\nu[n]]=0.
\nan

Given a finite dimensional  irreducible $Y(gl(M))\dot+Y(gl(N))$ module
$V_0$,  we define the action of the following $\Y$ generators
$\{t^i_j[n],\ t^\mu_\nu[n], \ t^i_\mu[n] |$ $i,   j=1, 2, ..., M$,
$\mu, \nu=M+1, M+2, ..., M+N$, $0<n\in{\bf Z}_+\}$ on $V_0$ by
\ban
t^i_\mu[n] v&=&0, \\
t^i_j[n] v &=& \phi^i_j[n] v, \\
t^\mu_\nu[n] v &=& \psi^\mu_\nu[n] v,  \ \  \forall v\in V_0.
\nan
It follows from the first equation and the $\Y$ defining relations
(\ref{modes}) that $[t^i_j[m],\ t^\mu_\nu[n]] v=0$,
$\forall v\in V_0$.  Thus this definition is self consistent.

We further define the vector space $\bar V$ spanned by
\ban
\left\{ \prod_{p\in\Phi_-^{(1)}}^\succ Q^{\{k_{p}\}}_{p}[\{n_{p}\}]\otimes v
\ | \ \forall \{k_{p}\}, \{n_{p}\}; \ \ v\in V_0\right\},
\nan
which furnishes an indecomposable $\Y$ module, with the module action
defined in the obvious way.  Then the  quotient module of $\bar V$ by its
unique maximal invariant submodule is the irreducible $\Y$ module which we
intend to construct.
{}From Theorems (\ref{uniqueness}) and (\ref{iff}) we deduce that
this construction yields all the  finite dimensional irreps of $\Y$.

Kac' construction\cite{Kac} proves to be useful for studying
the representation theory of Lie superalgebras.
A generalization of  the method also enabled us
to develop a relatively satisfactory representation theory for the
type I quantum supergroups\cite{Zhang}.  We hope that the construction
presented here will also provides a practicable method for
investigating detailed structures of the finite dimensional irreps of $\Y$.

\vspace{2cm}
\noindent
{\bf Acknowledgements}:  Part of this work was done  while I
visited the Institutes of Applied Mathematics and Theoretical Physics,
Chinese Academy of Sciences, Beijing.
I wish to thank Professors S. K. Wang and K. Wu for their hospitality.
Financial support from the Australian Research Council
is gratefully acknowledged.

\pagebreak

\end{document}